\documentclass[a4paper, 12pt]{amsart}
\usepackage[margin = 3.0cm]{geometry}
\usepackage{amsthm}
\usepackage[colorlinks = true, citecolor = blue, linkcolor = blue]{hyperref}
\usepackage{microtype}
\usepackage{booktabs}
\usepackage{tikz}

\usepackage[scr = euler]{mathalpha}
\newtheorem{thmcounter}{none}[section]
\theoremstyle{definition}
\newtheorem{definition}[thmcounter]{Definition}
\newtheorem{example}[thmcounter]{Example}
\newtheorem{convention}[thmcounter]{Convention}
\newtheorem{remark}[thmcounter]{Remark}
\theoremstyle{plain}
\newtheorem{lemma}[thmcounter]{Lemma}
\newtheorem{proposition}[thmcounter]{Proposition}
\newtheorem{theorem}[thmcounter]{Theorem}
\newtheorem{corollary}[thmcounter]{Corollary}

\newtheorem{observation}[thmcounter]{Observation}

\everymath{\displaystyle}
\counterwithin{equation}{section}
\counterwithin{table}{section}
\counterwithin{figure}{section}

\newcommand{\cat}{\operatorname{Cat}}

\newcommand{\su}{\operatorname{SU}}
\newcommand{\rank}{\operatorname{rank}}
\newcommand{\brank}{\operatorname{\underline{rank}}}

\newcommand{\diff}{\mathrm{D}}
\newcommand{\df}{\mathrm{d}}
\newcommand{\dist}{\mathrm{dist}}
\newcommand{\res}{\mathrm{Res}}
\newcommand{\disc}{\operatorname{Disc}}

\renewcommand{\Re}{\operatorname{Re}}
\renewcommand{\Im}{\operatorname{Im}}

\title{On the Real Waring Rank of a Biharmonic Form}
\author{Hua-lin Huang}
\address{School of Mathematical Sciences, Huaqiao University, Quanzhou 362021, China}
\email{hualin.huang@hqu.edu.cn}
\author{Yiming Liu}
\address{School of Mathematical Sciences, Huaqiao University, Quanzhou 362021, China}
\email{yimingliu@hqu.stu.edu.cn}
\author{Yilun Tang}
\address{School of Mathematical Sciences, Huaqiao University, Quanzhou 362021, China}
\email{18106088237@163.com}
\date{\today}

\begin{document}
    \begin{abstract}
        \par We develop an effective method for determining the real Waring rank and the real border rank of biharmonic binary forms. By means of real apolarity, we convert both of these problems into membership questions for regions bounded by a hypocycloid and an epicycloid. We further provide a computable algebraic criterion for the relevant memberships, which involves only the signs of certain real-valued resultants.
    \end{abstract}

    \maketitle

    \section{Introduction}

\par \noindent The Waring problem for polynomials is a classical topic in algebraic geometry and invariant theory that concerns the representation of forms as sums of powers of linear forms. Within this broad class of problems, the binary case is particularly tractable, and the complex rank of binary forms has been determined in \cite{Comas-2011, Iar-1999}. By contrast, substantially less is known about real rank, and most of the related problems remain open. In this paper we study the real rank and the real border rank of biharmonic binary forms, and we show that their rank stratification is governed by two cycloidal regions in the complex plane and that explicit algebraic criteria can be derived in terms of resultants.

\subsection{Problem statement} Let $f \in \mathbb{R}[x, y]_d$ be a real binary form of degree $d$. A representation of $f$ as
$$
    f = \lambda_1 (a_1 x + b_1 y)^d + \cdots + \lambda_r (a_r x + b_r y)^d,
$$
where $\lambda_j \in \mathbb{R}^{*}$ and the real linear forms $a_j x + b_j y$ are pairwise non-proportional, is called a \emph{Waring decomposition} of $f$, and the integer $r$ is called the \emph{length} of this decomposition. The \emph{real rank} of $f$, denoted by $\rank_\mathbb{R} f$, is defined to be the smallest integer $r$ such that $f$ admits a Waring decomposition of length $r$. A concept closely related to the real rank is the border rank. The real border rank of $f$ is defined as
$$
    \brank_\mathbb{R} f := \min \left\{
        s : f \in \overline{\left\{g : \rank_\mathbb{R} g \le s\right\}}
    \right\},
$$
where the closure is taken in the Euclidean topology on the coefficient space.

\par A real binary form $f \in \mathbb{R}[x, y]_d$ is said to be \emph{biharmonic} if $\Delta^2 f = 0$, where $\Delta = \partial_x^2 + \partial_y^2$ denotes the Laplace operator. The problem considered in this paper is to determine the real rank and the real border rank of every nonzero biharmonic binary form. In particular, we ask which values of these two invariants can occur, when they differ from each other, and whether they can be determined explicitly from the coefficients of the given form.

\par For convenience, we call a form \emph{real-rooted} if it splits completely into real linear factors. The main difficulty of this problem stems from real apolarity: a binary form $f$ admits a Waring decomposition of length $s$ if and only if the $s$-th homogeneous component $(f^\bot)_s$ of its apolar ideal $f^\bot$ contains a square-free real-rooted form. Similarly, the real border rank is governed by the existence of a real-rooted (not necessarily square-free) form in $(f^\bot)_s$. In general, an element of $(f^\bot)_s$ depends on several parameters, so determining whether some choice of these parameters yields a real-rooted or a square-free form is typically a multiparameter semialgebraic feasibility problem. The main challenge is therefore to exploit the special structure imposed by biharmonicity so as to render these real-rootedness and square-freeness conditions explicit.

\subsection{Related works} The study of representing binary forms as sums of powers of linear forms goes back to Sylvester and Gundelfinger \cite{Syl-1904A, Syl-1904B, Gund-1887}, in which the apolar approach to the Waring problem originated. For complex forms in any number of variables, the generic rank was determined by Alexander and Hirschowitz \cite{Alex-1995}. In the binary case, Comas and Seiguer described the complex rank strata and obtained an effective rank algorithm \cite{Comas-2011}. Huang, Miao and Ye \cite{Huang-2025A} gave an elementary treatment of the Waring problem for complex binary forms.

\par The behavior of the rank over the real field is substantially different from that over the complex field \cite{Rez-2013, Tok-2017, Rez-2017}. Typical real ranks were studied in low degrees by Comon and Ottaviani \cite{Comon-2012} and in higher degrees by Ballico \cite{Ball-2014}, before being completely determined by Blekherman \cite{Blek-2015}. The algebraic boundaries separating typical-rank regions were investigated by Brambilla and Staglianò \cite{Bram-2018, Bram-2021}. The relation between real-rootedness and maximal rank was studied by Causa and Re \cite{Cau-2011}, while Blekherman and Sinn \cite{Blek-2016} proved that a real binary form that is not a power of a linear form attains maximal real rank precisely when it is real-rooted. Recently, Huang, Tang, Ye, and Zhu \cite{Huang-2025B} proved that every harmonic binary form has Waring rank equal to its degree and gave a constructive decomposition procedure.

\subsection{Main approach and results} Our approach begins with real apolarity. Let $f \in \mathbb{R}[x, y]_d$ be a nonzero biharmonic form of degree $d$. Note that the biharmonic condition imposes no restriction when $d \le 3$; these cases are easy to handle because of the small degree. We therefore concentrate on degrees $d \ge 4$. By analyzing the homogeneous components of the apolar ideal $f^\bot$, we show that 
$$
    \rank_\mathbb{R} f, \ \brank_\mathbb{R} f \in \{d - 2, d - 1, d\}.
$$
Hence, it suffices to determine whether the $(d - 2)$-th and $(d - 1)$-th homogeneous components of $f^\bot$ contain a real-rooted (or additionally, square-free) form.

\par We mainly focus on the non-harmonic case, since in the harmonic case both the rank and the border rank attain their maximal value $d$. In this situation, the coefficients of $f$ determine a point $T_f$ in the complex plane $\mathbb{C}$. We prove that the required real-rootedness and square-freeness conditions in the relevant homogeneous components can be characterized equivalently by the position of $T_f$ relative to the $(d - 2)$-cusped hypocycloid $\Pi_{d - 2}$ and the $(d - 2)$-cusped epicycloid $\Gamma_{d - 2}$. Specifically, we have the following rank classification (see Theorem \ref{theorem:region-crit}):
$$
    \rank_\mathbb{R} f = \begin{cases}
        d - 2, & T_f \in \Omega_{d - 2}^{-}, \\
        d - 1, & T_f \in \Omega_{d - 2}^{+} \setminus \Omega_{d - 2}^{-}, \\
        d, & \text{otherwise}.
    \end{cases}
$$
and 
$$
    \brank_\mathbb{R} f = \begin{cases}
        d - 2, & T_f \in \Omega_{d - 2}^{-} \sqcup \Pi_{d - 2}, \\
        d - 1, & T_f \in (\Omega_{d - 2}^{+} \sqcup \Gamma_{d - 2}) \setminus (\Omega_{d - 2}^{-} \sqcup \Pi_{d - 2}), \\
        d, & \text{otherwise}.
    \end{cases}
$$
Here $\Omega_{d - 2}^-$ and $\Omega_{d - 2}^+$ denote the regions enclosed by $\Pi_{d - 2}$ and $\Gamma_{d - 2}$, respectively. For $d = 4$, we set $\Pi_2 = [-2, 2]$ and $\Omega_2^- = (-2, 2)$. This classification reduces the problem of determining the two ranks to the corresponding region-membership problems.

\par We further provide easily verifiable algebraic criteria for the above region-membership problems. Together with the classification, these criteria yield an explicit algorithm for computing the two ranks of biharmonic forms. Specifically, for the point $T_f$ above, we consider the univariate polynomials
$$
    P_{d - 2, T_f}(u) = u^{d - 1} - (d - 1)u + T_f
    \quad\text{and}\quad
    Q_{d - 2, T_f}(u) = (d - 3)u^{d - 2} - T_f u^{d - 3} + 1,
$$
and their conjugate reciprocal polynomials $P_{d - 2, T_f}^\#$ and $Q_{d - 2, T_f}^\#$, respectively. Then, the real rank and the real border rank of $f$ can be determined in the following way (see Theorem \ref{theorem:maintheorem-dge4}):
$$
    \rank_\mathbb{R} f = \begin{cases}
        d - 2, & \res(Q_{d - 2, T_f}, Q^\#_{d - 2, T_f}) > 0, \\
        d - 1, & \text{$(-1)^{d - 2}\res(P_{d - 2, T_f}, P^\#_{d - 2, T_f}) > 0$, $\res(Q_{d - 2, T_f}, Q^\#_{d - 2, T_f}) \le 0$}, \\
        d, & (-1)^{d - 2}\res(P_{d - 2, T_f}, P^\#_{d - 2, T_f}) \le 0,
    \end{cases}
$$
and
$$
    \brank_\mathbb{R} f = \begin{cases}
        d - 2, & \res(Q_{d - 2, T_f}, Q^\#_{d - 2, T_f}) \ge 0, \\
        d - 1, & \text{$(-1)^{d - 2}\res(P_{d - 2, T_f}, P^\#_{d - 2, T_f}) \ge 0$, $\res(Q_{d - 2, T_f}, Q^\#_{d - 2, T_f}) < 0$}, \\
        d, & (-1)^{d - 2}\res(P_{d - 2, T_f}, P^\#_{d - 2, T_f}) < 0.
    \end{cases}
$$
Here, when $d = 4$, the conditions $\res(Q_{2, T_f}, Q_{2, T_f}^{\#}) > 0$ and $\res(Q_{2, T_f}, Q_{2, T_f}^{\#}) \ge 0$ should be replaced by $T_f \in (-2, 2)$ and $T_f \in [-2, 2]$, respectively.

\par Finally, we deduce from the above classification that a biharmonic form of degree $d \ge 4$ has maximal real rank $d$ if and only if it is real-rooted, in agreement with the criterion of Blekherman and Sinn \cite[Theorem 2.2]{Blek-2016}. 

\subsection{Organization of this paper} In Section 2, we review some basic results on apolarity and biharmonic binary forms and treat the low-degree cases. Section 3 reduces the real rank and real border rank problems to membership in two cycloidal regions and determines their geometry. Section 4 provides effective algebraic criteria for the resulting classification, establishes the equivalence between real-rootedness and maximal real rank, and concludes with illustrative examples.
    \section{Apolarity and Biharmonic Forms} 

\par \noindent This section gathers some preliminaries on apolarity and biharmonic binary forms. Before proceeding, we need to fix the following convention:

\begin{convention}
    \par A nonzero real binary form $f \in \mathbb{R}[x, y]_d$ is said to be \emph{real-rooted} if it splits completely into real linear factors.
\end{convention}

\subsection{Apolarity for real binary forms} Let $\mathbb{K}$ be a field of characteristic zero, and consider the \emph{partial derivative action} of $\mathbb{K}[X, Y]$ on $\mathbb{K}[x, y]$ defined as
$$
    \varphi \circ f = \varphi(\partial_x, \partial_y) f, \quad \varphi \in \mathbb{K}[X, Y], ~ f \in \mathbb{K}[x, y].
$$
For any binary form $f \in \mathbb{K}[x, y]_d$, its \emph{apolar ideal} is defined to be
$$
    f^\bot := \left\{\varphi \in \mathbb{K}[X, Y] : \varphi \circ f = 0\right\}.
$$
One can immediately check that $f^\bot$ is a homogeneous ideal, and for any $s \le d$, its $s$-th homogeneous component is exactly the kernel of the $\mathbb{K}$-homomorphism
$$
    \cat_f(s) : \mathbb{K}[X, Y]_s \to \mathbb{K}[x, y]_{d - s}, \quad \varphi \mapsto \varphi \circ f.
$$
More precisely, writing
$$
    f(x, y) = \binom{d}{0} a_0 x^d + \binom{d}{1} a_1 x^{d - 1}y + \cdots + \binom{d}{d} a_d y^d,
$$
the form $\varphi(X, Y) = b_0 X^s + b_1 X^{s - 1}Y + \cdots + b_s Y^s$ belongs to $(f^\bot)_s$ if and only if its coefficient vector $(b_0, b_1, \cdots, b_s)$ is a solution of the linear system.
$$
    \left(\begin{array}{ccccc}
        a_0 & a_1 & \cdots & a_s \\
        a_1 & a_2 & \cdots & a_{s + 1} \\
        \vdots & \vdots & & \vdots \\
        a_{d - s} & a_{d - s + 1} & \cdots & a_d
    \end{array}\right)
    \left(\begin{array}{c}
        b_0 \\ b_1 \\ \vdots \\ b_s
    \end{array}\right) = 0.
$$

\begin{theorem}[Structure Theorem {\cite[Theorem 1.44]{Iar-1999}}]
    \label{theorem:struct-apid}
    \par Let $f \in \mathbb{K}[x, y]_d$. Then its apolar ideal $f^\bot$ is generated by two coprime forms $\psi_1$ and $\psi_2$ with
    $$
        \deg \psi_1 + \deg \psi_2 = d + 2.
    $$
    Moreover, if $\deg \psi_1 < \deg \psi_2$, then $\psi_1$ is uniquely determined up to scalar.
\end{theorem}

\begin{lemma}[Generalized Apolarity Lemma {\cite[Theorem 1.31]{Iar-1999}}]
    \label{lemma:gen-aplem}
    \par Let $f \in \mathbb{R}[x, y]_d$, and let $b_1 X - a_1 Y, \cdots, b_m X - a_m Y \in \mathbb{R}[X, Y]_1$ be pairwise non-proportional linear forms. Then the product 
    $$
        \prod_{k = 1}^{m} \left(b_k X - a_k Y\right)^{s_k} \quad \text{with} \quad s_1 + \cdots + s_m \le d
    $$
    belongs to $f^\bot$ if and only if there exist $g_k \in \mathbb{R}[x, y]_{s_k - 1}~(k = 1, \cdots, m)$ such that $f$ can be expressed as the sum
    $$
        f = \sum_{k = 1}^{m} g_k \cdot \left(a_k x + b_k y\right)^{d - s_k + 1}.
    $$
\end{lemma}

\begin{corollary}
    \label{corollary:homcpn-apid}
    \par Let $f \in \mathbb{R}[x, y]_d$. Then
    \begin{enumerate}
        \item $\rank_\mathbb{R} f \le s$ if and only if there exists a square-free real-rooted form in $(f^\bot)_s$.
        \item $\brank_\mathbb{R} f \le s$ if and only if there exists a real-rooted (not necessarily square-free) form in $(f^\bot)_s$.
    \end{enumerate}
\end{corollary}

\par A version of this result appears in \cite[Lemma 1.38]{Iar-1999}, where it is proved by algebro-geometric arguments over an algebraically closed field. Therefore, that proof is not well suited to the present situation. Here, we give instead a self-contained, elementary proof that avoids such machinery.

\begin{proof}[Proof of Corollary \ref{corollary:homcpn-apid}]
    \par Part (1) is the usual real apolarity criterion, and we only prove part (2). Throughout the proof, closure is taken in the Euclidean topology on the coefficient space of $\mathbb{R}[x, y]_d$.

    \par \textsc{The ``only if'' part}. Suppose that $\brank_\mathbb{R} f \le s$. Then there exists a sequence $\{f_n\} \subseteq \mathbb{R}[x, y]_d$ converging to $f$ such that $\rank_\mathbb{R} f_n \le s$ for every $n$. By part (1), after multiplying by several non-proportional real linear factors when necessary, we may choose a square-free real-rooted form $\varphi_n \in (f_n^\bot)_s$ for every $n$. Moreover, by rescaling $\varphi_n$, we may assume that their coefficient vectors lie on the real unit sphere
    $$
        \mathbb{S}^s(\mathbb{R}) := \left\{(x_0, \cdots, x_s) \in \mathbb{R}^{s + 1} : |x_0|^2 + \cdots + |x_s|^2 = 1\right\}.
    $$
    Since $\mathbb{S}^s(\mathbb{R})$ is compact in $\mathbb{R}^{s + 1}$, after passing to a subsequence if necessary, we may assume that $\varphi_{n}$ converges to $\varphi$. Note that the restriction of the action $\circ$ on $\mathbb{R}[X, Y]_s \times \mathbb{R}[x, y]_d$ is continuous. Thus, we have $\varphi \circ f = 0$, and it remains to show that $\varphi$ is a product of real linear forms.

    \par Indeed, assume to the contrary that $\varphi$ has an irreducible quadratic factor, and consider the dehomogenizations
    $$
        \tilde{\varphi}(x) = \varphi(x, 1) \quad \text{and} \quad
        \tilde{\varphi}_{n}(x) = \varphi_{n}(x, 1).
    $$
    Then, the sequence $\{\tilde{\varphi}_{n}\}$ converges to $\tilde{\varphi}$ and all $\tilde{\varphi}_{n}$ are real-rooted, while $\tilde{\varphi}(x)$ possesses at least one non-real root $z_0$. Now, consider the circle
    $$
        C_{\varepsilon} := \left\{z \in \mathbb{C} : |z - z_0| = \varepsilon\right\}.
    $$
    When $\varepsilon$ is sufficiently small, we may assume that $C_\varepsilon$ and its interior contain no roots of $\tilde{\varphi}$ except $z_0$, and that $C_\varepsilon$ does not intersect the real axis. Then, by Rouch\'e's Theorem \cite[Theorem 6.24]{Rod-2013}, for sufficiently large $n \in \mathbb{Z}_{>0}$, the polynomial $\tilde{\varphi}_{n}(x)$ has at least one zero inside $C_\varepsilon$. This implies that $\varphi_{n}(x, y)$ admits an irreducible quadratic factor, which is a contradiction.

    \par \textsc{The ``if'' part}. By Lemma \ref{lemma:gen-aplem}, we may suppose that $f$ admits a representation as follows:
    $$
        f = g_1 \cdot (a_1 x + b_1 y)^{d - s_1 + 1} + \cdots + g_m \cdot (a_m x + b_m y)^{d - s_m + 1},
    $$
    where $s = s_1 + \cdots + s_m$ and $a_k x + b_k y$ are pairwise non-proportional real linear forms. For each term $g_k \cdot (a_k x + b_k y)^{d - s_k + 1}$ in the representation above, after a suitable real linear change of coordinates, we may assume it has the form
    $$
        \gamma_k(u, v) \cdot u^{d - s_k + 1} = c_0 u^d + c_1 u^{d - 1} v + \cdots + c_{s_k - 1} u^{d - s_k + 1} v^{s_k - 1}.
    $$
    Expressing the higher-order derivatives in terms of finite differences, we see that for each $t \in \{1, \cdots, s_k\}$, we have
    $$
        \begin{aligned}
            u^{d - t + 1} v^{t - 1} &= \frac{\partial^{t - 1}}{\partial \varepsilon^{t - 1}} \left(\frac{(u + \varepsilon v)^d}{d! / (d - t + 1)!}\right) \bigg|_{\varepsilon = 0} \\ 
            &= \frac{(d - t + 1)!}{d!} \lim_{\varepsilon \to 0} \frac{1}{\varepsilon^{t - 1}} \sum_{j = 0}^{t - 1} (-1)^{t - 1 - j} \binom{t - 1}{j} (u + j\varepsilon v)^d.
        \end{aligned}
    $$
    Thus each term $\gamma_k(u, v) \cdot u^{d - s_k + 1}$ is the limit of $\mathbb{R}$-linear combinations of the $s_k$ powers
    $$
        u^d, \quad (u + \varepsilon_n v)^d, \quad \cdots, \quad (u + (s_k - 1)\varepsilon_n v)^d,
    $$
    where $\varepsilon_n \to 0$ as $n \to +\infty$. Now, to complete the proof, we only need to ensure that all these $d$-th powers are pairwisely non-proportional, and for this it suffices to suitably adjust the rate at which the corresponding $\varepsilon_n$ tends to zero.
\end{proof}

\subsection{Biharmonic forms} A real binary form $f \in \mathbb{R}[x, y]_d$ is said to be \emph{biharmonic} if it is annihilated by the square of the Laplace operator $\Delta$, that is, 
$$
    \Delta^2 f = \left(\partial_x^2 + \partial_y^2\right)^2 f = 0.
$$

\begin{observation}
    \label{observation:biharmform}
    \par After the complex change of variables $z = x + iy$ and $\bar{z} = x - iy$, any biharmonic binary form $f \in \mathbb{R}[x, y]_d$ can be expressed as
    $$
        f = \binom{d}{0} a z^d + \binom{d}{0} \bar{a} \bar{z}^d + \binom{d}{1} b z^{d - 1} \bar{z} + \binom{d}{d - 1} \bar{b} z \bar{z}^{d - 1}, 
        \quad a, b \in \mathbb{C}.
    $$
\end{observation}

\begin{proof}
    \par Note that any real binary form can be expressed as
    $$
        f = \binom{d}{0} c_0 z^d + \binom{d}{1} c_1 z^{d - 1} \bar{z} + \cdots + \binom{d}{d} c_d \bar{z}^d,
    $$
    where $c_j$ and $c_{d - j}$ are conjugate (i.e. $c_j = \bar{c}_{d - j}$). Since $\Delta = 4 \partial_z \partial_{\bar{z}}$, the biharmonic condition $\Delta^2 f = 0$ forces the corresponding coefficients $c_2, \cdots, c_{d - 2}$ to vanish.
\end{proof}

\par Now, let us discuss the possible real rank and real border rank of a biharmonic binary form. Let $f \in \mathbb{R}[x, y]_d$ be a nonzero biharmonic form of degree $d$. Then its real rank and real border rank are at most $d$ \cite[Theorem 4.10]{Rez-2013}. 

\par We first consider the case $d \le 3$, in which the biharmonic condition $\Delta^2 f = 0$ holds vacuously. Let $\disc(f)$ denote the discriminant of $f$.
\begin{itemize}
    \item \textsc{The case $d = 1$}. $f$ is a linear form and thus $\rank_\mathbb{R} f = \brank_\mathbb{R} f = 1$.
    \item \textsc{The case $d = 2$}. $f$ is a quadratic binary form and thus
    $$
        \rank_\mathbb{R} f = \brank_\mathbb{R} f = \begin{cases}
            1, & \disc(f) = 0, \\
            2, & \text{otherwise}.
        \end{cases}
    $$
    \item \textsc{The case $d = 3$}. This case has a concise description if one uses the results of Blekherman and Sinn \cite[Theorem 2.2]{Blek-2016}. More precisely, we have
    $$
        \rank_\mathbb{R} f = \begin{cases}
            1, & \text{$f$ is the cube of a real linear form}, \\
            2, & \disc(f) < 0, \\
            3, & \text{otherwise}.
        \end{cases}
    $$
    and
    $$
        \brank_\mathbb{R} f = \begin{cases}
            1, & \text{$f$ is the cube of a real linear form}, \\
            2, & \text{$\disc(f) \le 0$ but $f$ is not the cube of a real linear form}, \\
            3, & \text{otherwise}.
        \end{cases}
    $$
    We give a brief explanation for the border-rank criterion. Indeed, for a nonzero real binary cubic $f$, the real border rank equals $1$ precisely when $f$ is the cube of a real linear form. Moreover, a binary cubic with a double real root can be expressed as a limit of real rank-$2$ forms, whereas any binary cubic of rank at most $2$ has nonpositive discriminant.
\end{itemize}

\par For the cases $d = 4$ and $d = 5$, biharmonicity is nontrivial. In this case, neither the real rank nor the real border rank can be $1$, since no nonzero $d$-th power of a real linear form is biharmonic. There is also no quintic biharmonic form of real border rank at most $2$. Indeed, if $\brank_\mathbb{R} f = 2$, then, as border rank $1$ has already been excluded, Corollary \ref{corollary:homcpn-apid} and Theorem \ref{theorem:struct-apid} give
$$
    f^\bot = \langle\psi_1, \psi_2\rangle, \quad
    (\deg \psi_1, \deg \psi_2) = (2, 5),
$$
where $\psi_1$ is real-rooted. On the other hand, since $(X^2 + Y^2)^2 \in f^\bot$ and its degree is smaller than $\deg \psi_2$, it must be divisible by $\psi_1$. As $X^2 + Y^2$ is irreducible over $\mathbb{R}$, this forces $\psi_1$ to be a nonzero scalar multiple of $X^2 + Y^2$, contradicting the real-rootedness of $\psi_1$.

\begin{proposition}
    \label{proposition:possible-rk-brk}
    \par Let $f \in \mathbb{R}[x, y]_d$ be a biharmonic form of degree $d \ge 4$. Then
    $$
        \brank_\mathbb{R} f, \ \rank_\mathbb{R} f \in \{d - 2, d - 1, d\}.
    $$
\end{proposition}

\begin{proof}
    \par From the discussion above, the proposition holds for $d = 4$ and $d = 5$. We now prove the case $d \ge 6$. By Corollary \ref{corollary:homcpn-apid}, it suffices to show that there exists no real-rooted form in $(f^\bot)_k$ with $k < d - 2$. Theorem \ref{theorem:struct-apid} allows us to assume that $f^\bot = \langle \psi_1, \psi_2\rangle$ with $\deg \psi_1 \le \deg \psi_2$.

    \par \textsc{Case 1. $\deg \psi_1 = \deg \psi_2$}. In this case, we must have $d = 6$. Otherwise, $(f^\bot)_4$ would vanish, which contradicts $\Delta^2 f = 0$. By Theorem \ref{theorem:struct-apid}, $\deg \psi_1 = 4$ and thus all homogeneous components of degree less than $4$ vanish.

    \par \textsc{Case 2. $\deg \psi_1 < \deg \psi_2$}. Since $d \ge 6$ and $\deg \psi_1 + \deg \psi_2 = d + 2$, we have $\deg \psi_2 > 4$, which implies that the degree-$4$ element $(X^2 + Y^2)^2 \in f^\bot$ must be divisible by $\psi_1$. Thus, we obtain $\deg \psi_1 \le 4$ and $\deg \psi_2 \ge d - 2$. Moreover, $\psi_1$ is divisible by the irreducible quadratic $X^2 + Y^2$, so every element of $f^\bot$ of degree smaller than $d - 2$ is a multiple of $\psi_1$ and therefore cannot be real-rooted. Consequently, in both cases no homogeneous component $(f^\bot)_k$ with $k < d - 2$ contains a real-rooted form, which proves the proposition.
\end{proof}
    \section{Cycloidal Regions}

\par \noindent In this section we study the real rank and the real border rank of a biharmonic binary form of degree $d \ge 4$, writing $r = d - 2$ for notational convenience.

\subsection{Reduction of the rank problems} Let $f \in \mathbb{R}[x, y]_d$ be a biharmonic binary form. By apolarity, the real border rank of $f$ is determined by the lowest degree of the homogeneous component of $f^\bot$ containing a real-rooted polynomial, while the real rank of $f$ is determined by the lowest degree of the homogeneous component of $f^\bot$ containing a real-rooted and square-free polynomial. On the other hand, by the discussion in Section 2, for a biharmonic form of degree $d \ge 4$, only the homogeneous components $(f^\bot)_r$, $(f^\bot)_{r + 1}$ and $(f^\bot)_{r + 2}$ need to be considered in order to determine its rank and border rank.

\par We now characterize when a real-rooted (and square-free) form $\varphi \in \mathbb{R}[X, Y]_r$ belongs to $(f^\bot)_r$, and we do the same for the homogeneous component $(f^\bot)_{r + 1}$. Rewrite $f$ as the form in Observation \ref{observation:biharmform}, that is,
$$
    f = a z^d + \bar{a} \bar{z}^d + d b z^{d - 1} \bar{z} + d \bar{b} z \bar{z}^{d - 1}.
$$
If $b = 0$, then $\Delta f = 0$, and Theorem \ref{theorem:struct-apid} implies that neither $(f^\bot)_r$ nor $(f^\bot)_{r + 1}$ contains a real-rooted form. It therefore remains to consider only the case $b \neq 0$.

\par We first consider the $r$-th homogeneous component. Writing $Z = (X - iY)/2$ and $\bar{Z} = (X + iY)/2$, we may express any $\varphi \in \mathbb{R}[X, Y]_r$ in the complex coordinates $Z$ and $\bar{Z}$ as
$$
    \varphi = c_0 Z^r + c_1 Z^{r - 1}\bar{Z} + \cdots + c_r \bar{Z}^r, \quad c_j = \bar{c}_{r - j}.
$$
A direct calculation shows that the form $\varphi$ belongs to $(f^\bot)_r$ if and only if its coefficients satisfy the equations
$$
    a c_0 + b c_1 = 0 \quad\text{and}\quad b c_0 + \bar{b} c_r = 0.
$$
On the other hand, the form $\varphi$ is real-rooted (and square-free) if and only if $c_0 \neq 0$ and there exist (distinct) $\zeta_1, \cdots, \zeta_r \in \mathbb{C}$ with each $|\zeta_j| = 1$ such that 
$$
    \varphi = c_0 (Z + \zeta_1\bar{Z}) \cdots (Z + \zeta_r\bar{Z}).
$$
By Vieta's formulas, $\varphi$ satisfies both conditions above if and only if there exist $\zeta_1, \cdots, \zeta_r \in \mathbb{C}$, with $|\zeta_j| = 1$ for every $j$ and with the $\zeta_j$ pairwise distinct whenever $\varphi$ is required to be square-free, such that 
$$
    \sum_{j = 1}^{r} \zeta_j = - \frac{a}{b}
    \quad\text{and}\quad
    \prod_{j = 1}^{r} \zeta_j = - \frac{b}{\bar{b}}.
$$
\par Taking an arbitrary $r$-th root $\eta$ of $-b / \bar{b}$ and setting each $\zeta_j = \eta \omega_j$, we see that the latter condition is equivalent to the existence of (distinct) $\omega_1, \cdots, \omega_r \in \mathbb{C}$ such that 
$$
    |\omega_1| = \cdots = |\omega_r| = 1, \quad
    \sum_{j = 1}^{r} \omega_j = - \frac{a}{b\eta}, \quad
    \prod_{j = 1}^{r} \omega_j = 1.
$$

\par As for the $(r + 1)$-th homogeneous component, an argument analogous to the one above shows that a form 
$$
    \varphi = c_0 Z^{r + 1} + c_1 Z^r \bar{Z} + \cdots + c_{r + 1} \bar{Z}^{r + 1} \in \mathbb{R}[X, Y]_{r + 1}
$$
is real-rooted and square-free and belongs to $(f^\bot)_{r + 1}$ if and only if there exist (distinct) $\zeta_1, \cdots, \zeta_{r + 1} \in \mathbb{C}$ with each $|\zeta_j| = 1$ such that
$$
    a + b \sum_{j = 1}^{r + 1} \zeta_j + \bar{b} \prod_{j = 1}^{r + 1} \zeta_j = 0.
$$
Taking an arbitrary $r$-th root $\eta$ of $-b / \bar{b}$ and setting each $\zeta_j = \eta \omega_j$, then we see that the latter condition is equivalent to the existence of (distinct) $\omega_1, \cdots, \omega_{r + 1} \in \mathbb{C}$ such that 
$$
    |\omega_1| = \cdots = |\omega_{r + 1}| = 1
    \quad\text{and}\quad
    \sum_{j = 1}^{r + 1} \omega_j - \prod_{j = 1}^{r + 1} \omega_j = -\frac{a}{b\eta}.
$$

\begin{definition}
    \label{definition:regions}
    \par We define the following regions in the complex plane:
    $$
        \mathscr{E}_r := \left\{ \sum_{j = 1}^{r + 1} \omega_j - \prod_{j = 1}^{r + 1} \omega_j : |\omega_j| = 1 \right\}, \quad
        \mathscr{H}_r := \left\{ \sum_{j = 1}^{r} \omega_j : |\omega_j| = 1, \ \prod_{j = 1}^{r} \omega_j = 1 \right\}.
    $$
    Let $\mathscr{E}_r^{\dist} \subseteq \mathscr{E}_r$ and $\mathscr{H}_r^{\dist} \subseteq \mathscr{H}_r$ be the subsets obtained by requiring the $\omega_j$ to be pairwisely distinct.
\end{definition}

\begin{remark}
    \par The region $\mathscr{H}_r$ in Definition \ref{definition:regions} can be viewed as the set consisting of the traces of all $r \times r$ special unitary matrices, that is
    $$
        \mathscr{H}_r := \left\{\operatorname{trace}(M) : M \in \su(r)\right\}.
    $$
    This is related to the region studied by Kaiser \cite{Kai-2006}.
\end{remark}

\begin{observation}
    \label{observation:region-inclusion}
    \par Fix $r \ge 2$. Then the four regions in Definition \ref{definition:regions} satisfy the inclusions:
    $$
        \mathscr{H}_r \subseteq \mathscr{E}_r
        \quad\text{and}\quad
        \mathscr{H}_r^\dist \subseteq \mathscr{E}_r^\dist.
    $$
\end{observation}

\begin{proof}
    \par For any $\omega_1, \cdots, \omega_r$ with each $|\omega_j| = 1$ and $\omega_1 \cdots \omega_r = 1$, choosing $\omega_{r + 1}$ on the unit circle different from each $\omega_j~(j = 1, \cdots, r)$, we have
    $$
        \sum_{j = 1}^{r} \omega_j = \sum_{j = 1}^{r + 1} \omega_j - \prod_{j = 1}^{r + 1} \omega_j.
    $$
    This proves the two inclusions.
\end{proof}

\begin{theorem}
    \label{theorem:region-crit}
    \par Let $f \in \mathbb{R}[x, y]_d$ be a nonzero biharmonic form, written in the form of Observation \ref{observation:biharmform}. If $b = 0$, then $\brank_\mathbb{R} f = \rank_\mathbb{R} f = d$. Otherwise, by taking an arbitrary $r$-th root $\eta$ of $- b /\bar{b}$ and setting $T_f = - a / (b\eta)$, we obtain
    $$
        \rank_\mathbb{R} f = \begin{cases}
            d - 2, & T_f \in \mathscr{H}_r^{\dist}, \\
            d - 1, & T_f \in \mathscr{E}_r^\dist \setminus \mathscr{H}_r^\dist,\\
            d, & \text{otherwise}.
        \end{cases}
        \quad\text{and}\quad
        \brank_\mathbb{R} f = \begin{cases}
            d - 2, & T_f \in \mathscr{H}_r, \\
            d - 1, & T_f \in \mathscr{E}_r \setminus \mathscr{H}_r, \\
            d, & \text{otherwise}.
        \end{cases}
    $$
\end{theorem}

\begin{proof}
    \par Recall that $\rank_\mathbb{R} f, \ \brank_\mathbb{R} f \le d$. If $b = 0$, then $\Delta f = 0$ and every element of $f^\bot$ of degree less than $d$ must be divisible by $X^2 + Y^2$, which forces $\rank_\mathbb{R} f = \brank_\mathbb{R} f = d$. Suppose now that $b \neq 0$. The discussion at the beginning of this subsection shows that $(f^\bot)_r$ contains a real-rooted (resp. real-rooted and square-free) form if and only if $T_f \in \mathscr{H}_r$ (resp. $T_f \in \mathscr{H}_r^\dist$), and that $(f^\bot)_{r + 1}$ contains such a form if and only if $T_f \in \mathscr{E}_r$ (resp. $T_f \in \mathscr{E}_r^\dist$). Combining this with Proposition \ref{proposition:possible-rk-brk} completes the proof.
\end{proof}

\begin{remark}
    \par The value of $T_f$ in Theorem \ref{theorem:region-crit} depends on the choice of the $r$-th root $\eta$ of $- b / \bar{b}$, but the conclusion is independent of this choice. Indeed, different choices of $\eta$ merely rotate the point $T_f$ by $r$-th roots of unity, while all the regions defined in Definition \ref{definition:regions} are invariant under these rotations.
\end{remark}

\par However, the criterion supplied by Theorem \ref{theorem:region-crit} is not yet convenient to apply directly. We need a more detailed description of the regions, together with an easily verifiable algebraic criterion for the corresponding region-membership problem. This is what we do next.

\subsection{Geometry of the regions \texorpdfstring{$\mathscr{E}_r$ and $\mathscr{H}_r$}{Er and Hr}} We study the geometric properties of the regions $\mathscr{E}_r$ and $\mathscr{H}_r$ by means of their parametrizations. Let us consider the following parameter domains:
$$
    \mathscr{M}_r = \left\{
        (\omega_1, \cdots, \omega_{r + 1}) : |\omega_j| = 1
    \right\}
$$
and 
$$
    \mathscr{N}_r = \left\{
        (\omega_1, \cdots, \omega_r) : |\omega_j| = 1, \prod_{j = 1}^{r} \omega_j = 1
    \right\}.
$$
Both $\mathscr{M}_r$ and $\mathscr{N}_r$ are compact in the Euclidean topology, and the regions $\mathscr{E}_r$ and $\mathscr{H}_r$ are precisely the images of the continuous maps
$$
    \mu_r(\omega_1, \cdots, \omega_{r + 1}) = \sum_{j = 1}^{r + 1} \omega_j - \prod_{j = 1}^{r + 1} \omega_j
    \quad\text{and}\quad
    \nu_r(\omega_1, \cdots, \omega_r) = \sum_{j = 1}^{r} \omega_j.
$$
Hence $\mathscr{E}_r$ and $\mathscr{H}_r$ are compact as well.

\begin{observation}
    \label{observation:cycloid}
    \par Let $T \in \mathbb{C}$.
    \begin{enumerate}
        \item Suppose that $T \in \mathscr{E}_r$ and $T = \mu_r(\omega)$ for some $\omega \in \mathscr{M}_r$. Then, $T$ is a point on the $r$-cusped epicycloid
        $$
            \Gamma_r := \left\{
                (r + 1) e^{it} - e^{i(r + 1)t} : 0 \le t \le 2\pi
            \right\}
        $$
        if $\omega = (\omega_1, \cdots, \omega_{r + 1})$ belongs to the subset
        $$
            \mathscr{X}_r := \left\{
                (p, \cdots, p) : |p| = 1
            \right\} 
            \cup \left\{
                \sigma(p, q, \cdots, q) : |p| = |q| = 1, q^r = 1, \sigma \in \mathfrak{S}_{r + 1}
            \right\}.
        $$

        \item Suppose that $T \in \mathscr{H}_r$ and $T = \nu_r(\omega)$ for some $\omega \in \mathscr{N}_r$. Then, $T$ is a point on the $r$-cusped hypocycloid 
        $$
            \Pi_r := \left\{
                (r - 1)e^{it} + e^{-i(r - 1)t} : 0 \le t \le 2\pi
            \right\}
        $$
        if $\omega = (\omega_1, \cdots, \omega_{r})$ belongs to the subset
        $$
            \mathscr{Y}_r := \left\{
                \sigma(p, \cdots, p, p^{-(r - 1)}) : |p| = 1, \sigma \in \mathfrak{S}_{r}
            \right\}.
        $$
    \end{enumerate}
\end{observation}

\begin{proof}
    \par A direct computation.
\end{proof}

\begin{lemma}
    \label{lemma:interior-er}
    \par If $\omega \in \mathscr{M}_r \setminus \mathscr{X}_r$, then $\mu_r(\omega)$ is an interior point of $\mathscr{E}_r$ and there exists another $\omega' = (\omega'_1, \cdots, \omega'_{r + 1})$ with pairwise distinct entries such that $\mu_r(\omega) = \mu_r(\omega')$.
\end{lemma}

\begin{proof}
    \par We may suppose that $\omega = (e^{i \theta_1}, \cdots, e^{i \theta_{r + 1}})$ with each $\theta_j \in [0, 2\pi]$. Consider the map
    $$
        U_r(t_1, \cdots, t_{r + 1}) = \left(
            \Re \mu_r(e^{i t_1}, \cdots, e^{i t_{r + 1}}),
            \Im \mu_r(e^{i t_1}, \cdots, e^{i t_{r + 1}})
        \right).
    $$
    The derivative of $U_r$ at the point $t = \theta$ is
    $$
        \diff U_r(\theta) = \left(\begin{array}{rcr}
            - \sin \theta_1 + \sin \left(\sum_{j = 1}^{r + 1}\theta_j\right) & \cdots & - \sin \theta_{r + 1} + \sin \left(\sum_{j = 1}^{r + 1}\theta_j\right) \\
            \cos \theta_1 - \cos \left(\sum_{j = 1}^{r + 1}\theta_j\right) & \cdots & \cos \theta_{r + 1} - \cos \left(\sum_{j = 1}^{r + 1}\theta_j\right) \\
        \end{array}\right).
    $$
    This matrix must be nonzero; otherwise each $e^{i\theta_j}$ would equal $e^{i(\theta_1 + \cdots + \theta_{r + 1})}$, so that all entries of $\omega$ coincide with some $p$ with $|p| = 1$. The identity $p = p^{r + 1}$ then gives $p^r = 1$, whence $\omega \in \mathscr{X}_r$, contrary to the hypothesis.
    
    \par We first treat the case in which $\diff U_r(\theta)$ has full real rank $2$. In this case, the Implicit Function Theorem provides an open neighborhood of $U_r(\theta)$ contained in $U_r(\mathbb{R}^{r + 1}) \subseteq \mathscr{E}_r$. Furthermore, if $(j, k)$ is a pair of indices such that $\omega_j = \omega_k$, then the $j$-th and $k$-th columns of $\diff U_r(\theta)$ coincide. Hence, $\mathbf{e}_j - \mathbf{e}_k \in \ker \diff U_r(\theta)$ and 
    $$
        \mathscr{K}_{j, k} := \left\{v = (v_1, \cdots, v_{r + 1}) \in \ker \diff U_r(\theta) : v_j = v_k\right\}
    $$
    is a proper subspace of $\ker \diff U_r(\theta)$. Since a finite union of proper linear subspaces cannot cover $\ker \diff U_r(\theta)$, we may choose $v \in \ker \diff U_r(\theta)$ such that $v_j \neq v_k$ whenever $\omega_j = \omega_k$. Choose a curve $\delta(s)$ in the level set $\{t : U_r(t) = U_r(\theta)\}$ with $\delta(0) = \theta$ and tangent to $v$ at $s = 0$. Then, for any pair $(\omega_j, \omega_k)$ with $\omega_j = \omega_k$, we have 
    $$
        \delta_j(s) - \delta_k(s) = s(v_j - v_k) + o(s)
    $$
    with $v_j - v_k \neq 0$. Hence, for all sufficiently small $s \neq 0$, the numbers $\omega'_j = e^{i\delta_j(s)}$ are pairwise distinct, while $U_r(\delta(s)) = U_r(\theta)$, so that $\mu_r(\omega) = \mu_r(\omega')$.
    
    \par Now, we assume that $\rank_\mathbb{R} \diff U_r(\theta) = 1$. Then, all the points $(\sin \theta_j, \cos \theta_j)$ lie on a common line $L$ passing through the point
    $$
        A = \left(\sin \left(\sum_{j = 1}^{r + 1} \theta_j\right), \cos \left(\sum_{j = 1}^{r + 1} \theta_j\right)\right).
    $$
    Since $A$ lies on the unit circle, let $B$ be the other point at which the line $L$ meets the unit circle, and denote by $n_A$ and $n_B$ the numbers of these points that coincide with $A$ and with $B$, respectively.

    \par We first claim that $2 \le n_A \le r$. Indeed, the condition $\omega \notin \mathscr{X}_r$ excludes the cases $n_A = 0$ and $n_A = r + 1$. For $n_A = 1$, we can assume without loss of generality that 
    $$
        \theta_1 = \cdots = \theta_r = \alpha
        \quad\text{and}\quad
        \theta_{r + 1} = \beta
    $$
    for some $\alpha, \beta \in [0, 2\pi]$ with $\alpha \neq \beta$. Then, comparing arguments with $\omega_{r + 1}$, we obtain $e^{i\beta} = e^{i(r\alpha + \beta)}$, which implies that $e^{ir\alpha} = 1$. However, the point $\omega$ obtained in this way necessarily belongs to $\mathscr{X}_r$, which yields a contradiction. Based on this claim, we may assume that $\theta_1 = \theta_2 = \xi$ such that $A = (\sin \xi, \cos \xi)$. By comparing the arguments with $\omega_1$, we obtain $e^{i(\theta_2 + \cdots + \theta_{r + 1})} = 1$ and
    $$
        \mu_r(\omega) =  e^{i\xi} + \sum_{j = 2}^{r + 1} \omega_j - e^{i\xi}  \prod_{j = 2}^{r + 1} \omega_j = \sum_{j = 2}^{r + 1} \omega_j,
    $$
    which shows that the value of $\mu_r(\omega)$ is independent of $\omega_1$. Therefore, if we replace $\theta_1$ by a value $\vartheta_1$ for which $(\sin \vartheta_1, \cos \vartheta_1)$ is neither $A$ nor $B$, then 
    $$
        \diff U_r(\vartheta_1, \theta_2, \cdots, \theta_{r + 1})
    $$
    has real rank $2$ without changing $U_r(\theta)$, equivalently without changing $\mu_r(\theta)$. This reduces to the case treated at the beginning of the proof.
\end{proof}

\begin{lemma}
    \label{lemma:interior-hr}
    \par If $\omega \in \mathscr{N}_r \setminus \mathscr{Y}_r$, then $\nu_r(\omega)$ is an interior point of $\mathscr{H}_r$ and there exists another $\omega' = (\omega'_1, \cdots, \omega'_r)$ with pairwise distinct entries such that $\nu_r(\omega) = \nu_r(\omega')$.
\end{lemma}

\begin{proof}
    \par The proof of this lemma is analogous in spirit to that of Lemma \ref{lemma:interior-er}. Suppose that $\omega = (e^{i\theta_1}, \cdots, e^{i\theta_r})$ for some $\theta_j \in [0, 2\pi]$, and consider the map
    $$
        V_r(t_1, \cdots, t_r) = \left(
            \sum_{j = 1}^{r} t_j, \sum_{j = 1}^{r} \cos t_j, \sum_{j = 1}^{r} \sin t_j
        \right).
    $$
    Its derivative at the point $t = \theta$ is
    $$
        \diff V_r (\theta_1, \cdots, \theta_r) = \left(
            \begin{array}{cccc}
                1 & 1 & \cdots & 1 \\
                - \sin \theta_1 & - \sin \theta_2 & \cdots & -\sin \theta_r \\
                \cos \theta_1 & \cos \theta_2 & \cdots & \cos \theta_r
            \end{array}
        \right).
    $$
    Note that $\omega$ does not lie in $\mathscr{Y}_r$. Thus, among these $e^{i\theta_j}$, either at least three distinct values occur, or exactly two distinct values occur with each appearing at least twice.

    \par In the first case, $\diff V_r(\theta)$ attains maximal real rank $3$. Then, by an argument similar to the one at the beginning of the proof of Lemma \ref{lemma:interior-er}, we conclude that $\nu_r(\omega)$ is an interior point of $\mathscr{H}_r$. Moreover, a suitable small perturbation of $\omega$ yields an $\omega'$ with pairwise distinct entries such that $\nu_r(\omega') = \nu_r(\omega)$.

    \par We now turn to the second case, in which we may assume that the entries of $\omega$ take precisely the two distinct values $e^{i(\theta + \alpha)}$ and $e^{i(\theta - \alpha)}$, with $\alpha \in (0, \pi)$; in particular, $\sin \alpha \neq 0$. Let us consider the function
    $$
        G(x, y) = \cos x + \cos y - 2 \cos \alpha.
    $$
    Note that $G(\alpha, \alpha) = 0$ and $G_y(\alpha, \alpha) = - \sin \alpha \neq 0$. Then the Implicit Function Theorem shows that for $x \neq \alpha$ sufficiently close to $\alpha$, there exists $y$ close to $\alpha$ such that $G(x, y) = G(\alpha, \alpha) = 0$ with $x, y \in (0, \pi)$ and $x \neq \pm y$. Consequently, we obtain four distinct values
    $$
        e^{i(\theta + x)}, \quad e^{i(\theta - x)}, \quad
        e^{i(\theta + y)}, \quad e^{i(\theta - y)},
    $$
    and their sum and product are
    $$
        2e^{i\theta}(\cos x + \cos y) = 4e^{i\theta} \cos \alpha
        \quad\text{and}\quad 
        e^{4i\theta},
    $$
    respectively, which coincide with those of two copies of $e^{i(\theta + \alpha)}$ and two copies of $e^{i(\theta - \alpha)}$. The problem is therefore reduced to the first case. We are done.
\end{proof}

\par For every $r \ge 2$, the epicycloid $\Gamma_r$ is a Jordan curve enclosing the origin. For $r \ge 3$, the same is true for the hypocycloid $\Pi_r$. By the Jordan Curve Theorem, we obtain the partitions
$$
    \mathbb{C} = \Omega_r^{+} \sqcup \Gamma_r \sqcup \mathrm{E}_r^{+}
    \quad\text{and}\quad
    \mathbb{C} = \Omega_r^{-} \sqcup \Pi_r \sqcup \mathrm{E}_r^{-}.
$$
Here $\Omega_r^-$ and $\Omega_r^+$ denote the bounded components, while $\mathrm{E}_r^-$ and $\mathrm{E}_r^+$ denote the unbounded ones. Based on the preceding discussion, we can now give a complete description of the regions $\mathscr{E}_r$ and $\mathscr{H}_r$.

\begin{theorem}
    \label{theorem:region-er-hr}
    \par The region $\mathscr{E}_r$ admits the description
    $$
        \mathscr{E}_r = \Omega_r^+ \sqcup \Gamma_r 
        \quad\text{with}\quad
        \mathscr{E}_r^\dist = \Omega_r^+.
    $$
    In addition, if $r \ge 3$, the region $\mathscr{H}_r$ admits the analogous description
    $$
        \mathscr{H}_r = \Omega_r^- \sqcup \Pi_r
        \quad\text{with}\quad
        \mathscr{H}_r^\dist = \Omega_r^-.
    $$
\end{theorem}

\begin{center}
    \begin{figure}[htbp]
        \label{figure:cycloid-region}
        \begin{tikzpicture}[scale = 0.5]
    % Epicycloid E_6
    \path[fill = gray!25, draw = gray!, line width = 1.0pt]
    plot[domain = 0:360, samples = 361, variable = \t]
    ({7*cos(\t) - cos(7*\t)}, {7*sin(\t) - sin(7*\t)}) -- cycle;

    % Hypocycloid H_6
    \path[fill = gray!50, draw = darkgray, line width = 1.0pt]
    plot[domain = 0:360, samples=361, variable=\t]
    ({5*cos(\t) + cos(5*\t)}, {5*sin(\t) - sin(5*\t)}) -- cycle;

    % Coordinate axes
    \draw[->, line width = 1.0pt] (-9,0) -- (9,0) node[right] {$\Re T$};
    \draw[->, line width = 1.0pt] (0,-9) -- (0,9) node[above] {$\Im T$};

    \foreach \x in {-8, -6, -4, -2, 2, 4, 6, 8}
    \draw (\x, 0.12) -- (\x, -0.12) node[below, font = \scriptsize] {$\x$};
    \foreach \y in {-8, -6, -4, -2, 2, 4, 6, 8}
    \draw (0.12, \y) -- (-0.12, \y) node[left, font = \scriptsize] {$\y$};

    % Region labels
    \node[inner sep = 1.5pt] at (1.2, 0.8) {$\Omega_6^{-}$};
    \node[inner sep = 1.5pt] at (1.2, 6.0) {$\Omega_6^{+}$};
    \node[inner sep = 1.5pt] at (5.4, 5.4) {$\Gamma_6$};
    \node[inner sep = 1.5pt] at (3.6, 2.4) {$\Pi_6$};
\end{tikzpicture}
        \caption{The regions enclosed by $\Gamma_6$ and $\Pi_6$.}
    \end{figure}
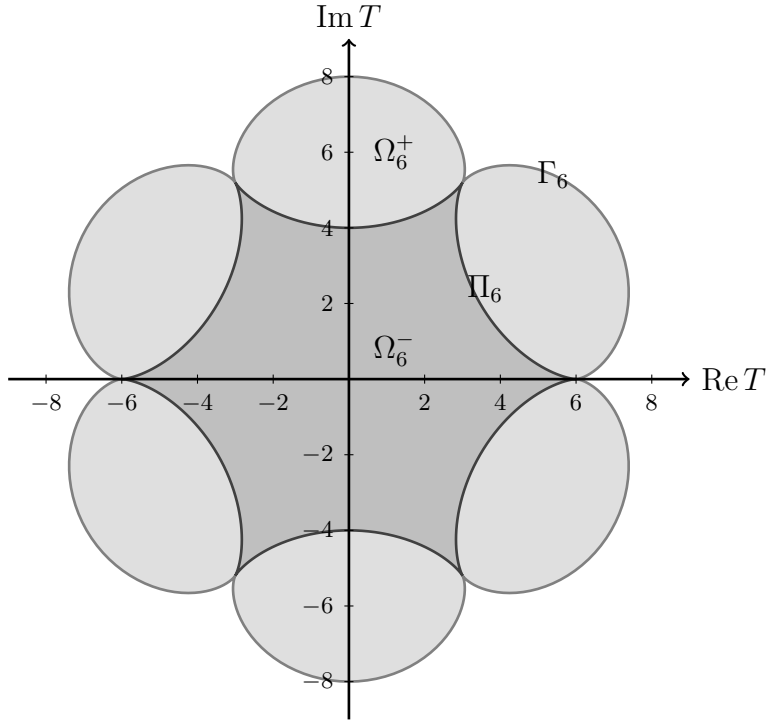
\end{center}

\begin{proof}
    \par We prove the theorem only for the region $\mathscr{E}_r$; the proof for $\mathscr{H}_r$ is analogous. 

    \par \textsc{Proof of $\mathscr{E}_r = \Omega_r^+ \sqcup \Gamma_r$}. Since $\mathscr{E}_r$ is compact, it suffices to show that $\partial \mathscr{E}_r \subseteq \Gamma_r \subseteq \mathscr{E}_r$ and that $\mathscr{E}_r \cap \Omega_r^+$ is nonempty. Let $T \in \partial \mathscr{E}_r$ and write $T = \mu_r(\omega)$ for some $\omega \in \mathscr{M}_r$. Being compact, $\mathscr{E}_r$ is the disjoint union of its interior and its boundary. Hence, by Lemma \ref{lemma:interior-er}, the corresponding $\omega$ must lie in $\mathscr{X}_r$, so that $T = \mu_r(\omega)$ belongs to the epicycloid $\Gamma_r$ by Observation \ref{observation:cycloid}. Conversely, the parametrization of the epicycloid shows immediately that every point of $\Gamma_r$ lies in $\mathscr{E}_r$, which proves the inclusions $\partial \mathscr{E}_r \subseteq \Gamma_r \subseteq \mathscr{E}_r$. On the other hand, taking $\omega_1, \cdots, \omega_r$ to be the roots of $z^r - (-1)^{r + 1}$ and $\omega_{r + 1} = 1$, we have
    $$
        \sum_{j = 1}^{r + 1} \omega_{j} - \prod_{j = 1}^{r + 1} \omega_{j} = 1 - 1 = 0.
    $$
    Hence $\mathscr{E}_r \cap \Omega_r^+$ contains the origin.

    \par \textsc{Proof of $\mathscr{E}_r^\dist = \Omega_r^+$}. Lemma \ref{lemma:interior-er} immediately gives the inclusion $\mathscr{E}_r^\dist \subseteq \Omega_r^+$. Conversely, let $T \in \Omega_r^+$ and write $T = \mu_r(\omega)$. By Observation \ref{observation:cycloid}, $\omega$ cannot lie in $\mathscr{X}_r$, so Lemma \ref{lemma:interior-er} provides another $\omega'$ with pairwise distinct entries such that $T = \mu_r(\omega')$. Thus $T \in \mathscr{E}_r^\dist$.
\end{proof}

\begin{remark}
    \par Theorem \ref{theorem:region-er-hr} does not apply to the region $\mathscr{H}_2$, since for $r = 2$, the hypocycloid $\Pi_2$ degenerates into the interval $[-2, 2]$ of the real axis. However, the description of $\mathscr{H}_2$ is straightforward, since $\mathscr{H}_2^{\dist}$ is precisely the open interval $(-2, 2)$.
\end{remark}

    \section{Criteria for Real Rank and Real Border Rank}

\par \noindent In Theorem \ref{theorem:region-crit}, we obtained a region-membership criterion for the real rank and the real border rank. Our goal in this section is to translate it into an easily computable algebraic criterion. Throughout this section, we keep the standing assumptions of Section 3, namely $d \ge 4$ and $r = d - 2$.

\subsection{Resultant criteria for real rank and real border rank} Let $f \in \mathbb{R}[x, y]_d$ be a nonzero biharmonic binary form of degree $d \ge 4$. Our aim in this subsection is to recast the region criterion in Theorem \ref{theorem:region-crit} as an algebraic criterion that is relatively straightforward to compute. Let $T$ be an arbitrary point in the complex plane. Consider the following polynomials
$$
    P_{r, T}(u) := u^{r + 1} - (r + 1)u + T
    \quad\text{and}\quad
    Q_{r, T}(u) := (r - 1)u^{r} - Tu^{r - 1} + 1.
$$
A direct substitution shows that $T$ lies on the epicycloid $\Gamma_r$ (resp. the hypocycloid $\Pi_r$) if and only if $P_{r, T}$ (resp. $Q_{r, T}$) has a zero on the unit circle $\partial \mathbb{D} = \left\{z \in \mathbb{C} : |z| = 1\right\}$.

\begin{remark}
    \label{remark:crit-windnum}
    \par For $T \in \mathbb{C}$ off the corresponding boundary cycloid, membership in $\Omega^+$ or $\Omega^-$ can be characterized in terms of the zero distribution of the polynomials $P_{r, T}$ and $Q_{r, T}$. Indeed, the regions $\Omega^+$ and $\Omega^-$ are enclosed by the Jordan curves $\Gamma_r$ and $\Pi_r$, respectively, and the winding number shows that $T$ lies in $\Omega^+$ if and only if
    $$
        \frac{1}{2\pi i} \oint_{\Gamma_r} \frac{1}{z - T} \df z
        = \frac{1}{2\pi i} \oint_{|u|=1} \frac{P'_{r, T}(u)}{P_{r, T}(u)} \df u = 1,
    $$
    which is equivalent to $P_{r, T}$ having exactly one zero, counted with multiplicity, in the unit disk $\mathbb{D} := \{z \in \mathbb{C} : |z| < 1\}$; while $T$ lies in $\Omega^-$ if and only if
    $$
        \frac{1}{2\pi i} \oint_{\Pi_r} \frac{1}{z - T} \df z
        = \frac{1}{2\pi i} \oint_{|u|=1} \frac{Q'_{r, T}(u)}{Q_{r, T}(u)} \df u - (r - 1) = 1,
    $$
    which is equivalent to all zeros of $Q_{r, T}$ lying in $\mathbb{D}$.
\end{remark}

\par Given a univariate polynomial $F = a_m u^m + \cdots + a_1 u + a_0 \in \mathbb{C}[u]$, we define its \emph{conjugate reciprocal polynomial}, denoted by $F^\#$, as
$$
    F^\# = \bar{a}_0 u^m + \cdots + \bar{a}_{m - 1} u + \bar{a}_m.
$$
We record the following basic properties of the conjugate reciprocal polynomial, which will be used below.

\begin{lemma}
    \label{lemma:properties-crp}
    \par Let $F = a_t u^t + \cdots + a_1 u + a_0 \in \mathbb{C}[u]$ with $a_0 a_t \neq 0$, and let $F^\#$ be its conjugate reciprocal polynomial.
    \begin{enumerate}
        \item A nonzero complex number $\alpha$ is a zero of $F$ of multiplicity $m$ if and only if $1 / \bar{\alpha}$ is a zero of $F^\#$ of the same multiplicity $m$.
        \item If $a_t a_0 \neq 0$, then the resultant $\res(F, F^\#)$ is real.
    \end{enumerate}
\end{lemma}

\begin{proof}
    \par Part (1) follows directly from the identity
    $$
        F^\#(u) = u^t \overline{F(1 / \bar{u})}.
    $$
    As for part (2), the resultant can be written as
    $$
        \res(F, F^\#) = a_t^t \prod_{j = 1}^{t} F^\#(\alpha_j)
        = |a_t|^{2t} \prod_{j = 1}^t (1 - |\alpha_j|^2) \prod_{1 \le j < k \le t} |1 - \alpha_j\bar{\alpha}_k|^2,
    $$
    which is evidently real.
\end{proof}

\par We now consider the conjugate reciprocal polynomials $P^\#_{r, T}$ and $Q^\#_{r, T}$ of $P_{r, T}$ and $Q_{r, T}$, respectively. By definition, they are given by
$$
    P_{r, T}^\# = \bar{T} u^{r + 1} - (r + 1) u^r + 1
    \quad\text{and}\quad
    Q_{r, T}^\# = u^r  - \bar{T} u + (r - 1).
$$

\begin{lemma}
    \label{lemma:resultant-crit}
    \par For every $T \in \mathbb{C}$, the two resultants $\res(P_{r, T}, P^\#_{r, T})$ and $\res(Q_{r, T}, Q^\#_{r, T})$ are real numbers. Moreover, for $r \ge 2$, we have
    $$
        (-1)^r \res(P_{r, T}, P^\#_{r, T})\begin{cases}
            > 0, & T \in \Omega_r^+, \\
            = 0, & T \in \Gamma_r, \\
            < 0, & T \in \mathrm{E}_r^+
        \end{cases}
    $$
    In addition, if $r \ge 3$, then
    $$
        \res(Q_{r, T}, Q^\#_{r, T})\begin{cases}
            > 0, & T \in \Omega_r^-, \\
            = 0, & T \in \Pi_r, \\
            < 0, & T \in \mathrm{E}_r^-.
        \end{cases}
    $$
\end{lemma}

\begin{proof}
    \par By Lemma \ref{lemma:properties-crp}(2), $\res(P_{r, T}, P^\#_{r, T})$ is real for all nonzero $T$ and $\res(Q_{r, T}, Q^\#_{r, T})$ is real for all $T$. Since $\res(P_{r, T}, P^\#_{r, T})$ depends continuously on $T$, passing to the limit shows that $\res(P_{r, 0}, P^\#_{r, 0})$ is real as well.

    \par We claim that each resultant vanishes if and only if $T$ lies on the corresponding cycloid. Indeed, if $T \in \Gamma_r$, say $T = (r + 1)e^{i\theta} - e^{i(r + 1)\theta}$ with $\theta \in [0, 2\pi]$, then $e^{i\theta}$ is a common zero of $P_{r, T}$ and $P^\#_{r, T}$. Conversely, write a common zero as $\rho e^{i\theta}$ with $\rho > 0$. By Lemma \ref{lemma:properties-crp}, both $\rho e^{i\theta}$ and $\rho^{-1} e^{i\theta}$ are zeros of $P_{r, T}$; eliminating $T$ gives
    $$
        (r + 1)(\rho - \rho^{-1}) = (\rho^{r + 1} - \rho^{- (r + 1)}) e^{ir\theta}.
    $$
    If $\rho \neq 1$, taking absolute values and dividing by $|\rho - \rho^{-1}|$ yields
    $$
        \rho^{r} + \rho^{r - 2} + \cdots + \rho^{-r} = (r + 1),
    $$
    while the arithmetic–geometric mean inequality shows that equality can hold only when $\rho = 1$, a contradiction. Hence $\rho = 1$. For the hypocycloid $\Pi_r$, the same argument gives
    $$
        (r - 1) (\rho - \rho^{-1}) e^{ir\theta} = \rho^{r - 1} - \rho^{-(r - 1)},
    $$
    and the arithmetic–geometric mean inequality again forces $\rho = 1$ for $r \ge 3$.

    \par Finally, expanding these two resultants as polynomials in $T$ and $\bar{T}$, we obtain
    $$
        \begin{cases}
            \res(P_{r, T}, P_{r, T}^\#) = (-1)^{r + 1} |T|^{2r + 2} + o(|T|^{2r + 2}), \\
            \res(Q_{r, T}, Q_{r, T}^\#) = - |T|^{2r - 2} + o(|T|^{2r - 2}).
        \end{cases}
    $$
    Hence both $(-1)^r\res(P_{r, T}, P_{r, T}^\#)$ and $\res(Q_{r, T}, Q_{r, T}^\#)$ are negative when $|T|$ is sufficiently large. On the other hand, a direct computation gives
    $$
        (-1)^r \cdot \res(P_{r, 0}, P_{r, 0}^\#) = (r^2 + 2r)^r > 0
    $$
    and, for $r \ge 3$,
    $$
        \res(Q_{r, 0}, Q_{r, 0}^\#) = (r^2 - 2r)^r > 0.
    $$
    Since each resultant is a continuous function of $T$ that does not vanish on any connected component of the complement of the corresponding cycloid, its sign is constant on each such component.
\end{proof}

\par With the preceding discussion, we are now in a position to present a criterion for the rank and border rank of biharmonic binary forms of degree $d \ge 4$.

\begin{theorem}
    \label{theorem:maintheorem-dge4}
    \par Let $f \in \mathbb{R}[x, y]_d$ be a nonzero biharmonic binary form of degree $d \ge 4$, written in the form given in Observation \ref{observation:biharmform}. If $b = 0$, then 
    $$
        \rank_\mathbb{R} f = \brank_\mathbb{R} f = d.
    $$
    Otherwise, with $T_f$ as in Theorem \ref{theorem:region-crit}, the following formulas hold:
    \begin{enumerate}
        \item For $d = 4$, we have
        $$
            \rank_\mathbb{R} f = \begin{cases}
                2, & T_f \in (-2, 2), \\
                3, & \text{$\res(P_{2, T_f}, P^\#_{2, T_f}) > 0$ and $T_f \notin (-2, 2)$}, \\
                4, & \res(P_{2, T_f}, P^\#_{2, T_f}) \le 0,
            \end{cases}
        $$
        and
        $$
            \brank_\mathbb{R} f = \begin{cases}
                2, & T_f \in [-2, 2], \\
                3, & \text{$\res(P_{2, T_f}, P^\#_{2, T_f}) \ge 0$ and $T_f \notin [-2, 2]$}, \\
                4, & \res(P_{2, T_f}, P^\#_{2, T_f}) < 0.
            \end{cases}
        $$
        \item For $d \ge 5$, we have
        $$
            \rank_\mathbb{R} f = \begin{cases}
                d - 2, & \res(Q_{r, T_f}, Q^\#_{r, T_f}) > 0, \\
                d - 1, & \text{$(-1)^r\res(P_{r, T_f}, P^\#_{r, T_f}) > 0$ and $\res(Q_{r, T_f}, Q^\#_{r, T_f}) \le 0$}, \\
                d, & (-1)^r\res(P_{r, T_f}, P^\#_{r, T_f}) \le 0.
            \end{cases}
        $$
        and
        $$
            \brank_\mathbb{R} f = \begin{cases}
                d - 2, & \res(Q_{r, T_f}, Q^\#_{r, T_f}) \ge 0, \\
                d - 1, & \text{$(-1)^r\res(P_{r, T_f}, P^\#_{r, T_f}) \ge 0$ and $\res(Q_{r, T_f}, Q^\#_{r, T_f}) < 0$}, \\
                d, & (-1)^r\res(P_{r, T_f}, P^\#_{r, T_f}) < 0.
            \end{cases}
        $$
    \end{enumerate}
\end{theorem}

\begin{proof}
    \par The result is an immediate consequence of Theorems \ref{theorem:region-crit} and \ref{theorem:region-er-hr} together with Lemma \ref{lemma:resultant-crit}.
\end{proof}

\begin{corollary}
    \label{corollary:real-rooted-biharm}
    \par Let $f \in \mathbb{R}[x, y]_d$ be a nonzero biharmonic binary form of degree $d \ge 4$. Then, in the notation of Observation \ref{observation:biharmform}, $f$ is real-rooted if and only if one of the following conditions holds:
    \begin{enumerate}
        \item[(i)] $b = 0$; 
        \item[(ii)] $b \neq 0$ and $(-1)^r \cdot \res(P_{r, T_f}, P^\#_{r, T_f}) \le 0$.
    \end{enumerate}
\end{corollary}

\begin{proof}
    \par Suppose first that $b = 0$. Since $f$ is nonzero, we have $a \neq 0$ and
    $$
        f(\cos \theta, \sin \theta) = 2 \Re (a e^{id\theta})
    $$
    If we treat this as a function of $\theta$, then it has exactly $d$ real projective zeros. Thus $f$ is a real-rooted form.
    
    \par Now suppose that $b \neq 0$. If $a = 0$, then $T_f$ is the origin and obviously lies in $\Omega_r^+$, so Lemma \ref{lemma:resultant-crit} gives $(-1)^r \cdot \res(P_{r, T_f}, P^\#_{r, T_f}) > 0$. On the other hand, note that
    $$
        f(\cos \theta, \sin \theta) = 2d \Re(be^{ir\theta})
    $$
    has exactly $r = d - 2$ real projective zeros. Hence neither condition (ii) nor real-rootedness holds in this case. It remains to consider the case $ab \neq 0$. Define the polynomial
    $$
        H_f(u) = a u^d + d b u^{d - 1} + d \bar{b} u + \bar{a}.
    $$
    Note that $H_f(e^{2i\theta}) = e^{id\theta} f(\cos \theta, \sin \theta)$. Therefore, $f$ is real-rooted if and only if all the zeros of $H_f$ lie on the unit circle. Since $H^\#_f = H_f$, Cohn's theorem (see \cite[pp. 120--121]{Cohn-1922} or \cite[Theorem A]{Cho-2013}) shows that it is equivalent to that all the zeros of $H'_f$ lie in the closed unit disk $\bar{\mathbb{D}} = \left\{z \in \mathbb{C} : |z| \le 1\right\}$. Choose $\eta$ such that $\eta^r = - b / \bar{b}$ and set $T_f = -a / (b\eta)$. A direct calculation gives 
    $$
        H'_f(u) = -db\eta u^{r + 1} P_{r, T_f}\left(\frac{1}{\eta u}\right).
    $$
    Since $|\eta| = 1$, all zeros of $H'_f$ lie in the closed unit disk $\bar{\mathbb{D}}$ if and only if $P_{r, T_f}$ has no zeros in the open unit disk $\mathbb{D}$. By Remark \ref{remark:crit-windnum}, the latter is equivalent to 
    $$
        T_f \in \Gamma_r \sqcup \mathrm{E}_r^+.
    $$
    This can be identified with $(-1)^r \cdot \res(P_{r, T_f}, P^\#_{r, T_f}) \le 0$ by Lemma \ref{lemma:resultant-crit}. 
\end{proof}

\begin{remark}
    \label{remark:blek-sinn}
    \par As observed in Section 2.2, no nonzero biharmonic form of degree $d \ge 4$ has real rank one. Hence the theorem of Blekherman and Sinn \cite[Theorem 2.2]{Blek-2016} applies and asserts that such a form has real rank $d$ if and only if it is real-rooted. This agrees with Corollary \ref{corollary:real-rooted-biharm}.
\end{remark}

\subsection{Some examples} We now illustrate the algebraic criteria established above by means of explicit examples of biharmonic binary forms. 

\begin{example}
    \par Consider the biharmonic forms $f \in \mathbb{R}[x, y]_d$ of degree $d \ge 6$ with $a = 0$ and $b \neq 0$ in the representation of Observation \ref{observation:biharmform}; such a form reads
    $$
        f = d b z^{d - 1} \bar{z} + d \bar{b} z \bar{z}^{d - 1}.
    $$
    For these forms the point $T_f = - a / (b\eta)$ is precisely the origin, which lies in the interior of $\Omega^-$; by Theorem \ref{theorem:region-crit}, their real rank is therefore $d - 2$.
\end{example}

\begin{example}
    \label{example:biharmform-1}
    \par Consider the biharmonic form
    $$
        f(x, y) = -140 x^6 y + 420 x^4 y^3 - 84 x^2 y^5 - 4 y^7.
    $$
    In the complex coordinates $z = x + iy$, $\bar{z} = x - iy$, and in the notation of Observation \ref{observation:biharmform}, this form becomes
    $$
        f = 5i z^7 - 5i \bar{z}^7 + 7i z^6 \bar{z} - 7i z \bar{z}^6.
    $$
    We may take $\eta = 1$, so that $T_f = -5$. Then
    $$
        (-1)^5 \cdot \res(P_{5, -5}, P_{5, -5}^\#) = 583 200 000 > 0
        \quad\text{and}\quad
        \res(Q_{5, -5}, Q_{5, -5}^\#) = -3 200 000 < 0,
    $$
    and Theorem \ref{theorem:maintheorem-dge4} yields $\rank_\mathbb{R} f = \brank_\mathbb{R} f = 6$.
    \begin{center}
        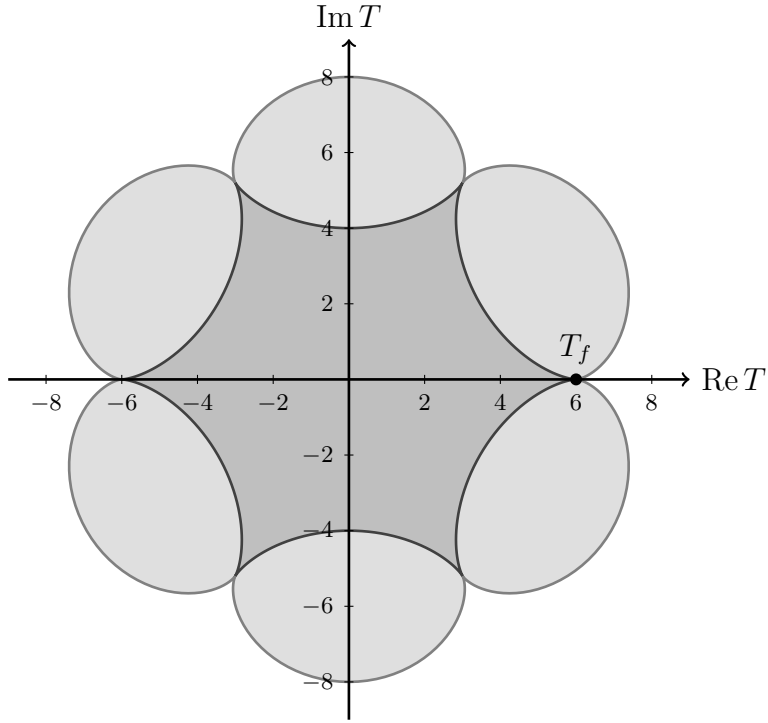
\begin{figure}[htbp]
            \label{figure:example-tf-1}
            \begin{tikzpicture}[scale = 0.5]
    % Epicycloid E_6
    \path[fill = gray!25, draw = gray, line width = 1.0pt]
    plot[domain = 0:360, samples = 361, variable = \t]
    ({6*cos(\t) - cos(6*\t)}, {6*sin(\t) - sin(6*\t)}) -- cycle;

    % Hypocycloid H_6
    \path[fill = gray!50, draw = darkgray, line width = 1.0pt]
    plot[domain = 0:360, samples=361, variable=\t]
    ({4*cos(\t) + cos(4*\t)}, {4*sin(\t) - sin(4*\t)}) -- cycle;

    % Coordinate axes
    \draw[->, line width = 1.0pt] (-8,0) -- (8,0) node[right] {$\Re T$};
    \draw[->, line width = 1.0pt] (0,-8) -- (0,8) node[above] {$\Im T$};

    \foreach \x in {-7, -5, -3, -1, 1, 3, 5, 7}
    \draw (\x, 0.12) -- (\x, -0.12) node[below, font = \scriptsize] {$\x$};
    \foreach \y in {-7, -5, -3, -1, 1, 3, 5, 7}
    \draw (0.12, \y) -- (-0.12, \y) node[left, font = \scriptsize] {$\y$};

    % Tf
    \fill (-5, 0) circle (4.5pt);
    \node[inner sep = 1.5pt] at (-5, 0.8) {$T_f$};
\end{tikzpicture}
            \caption{The location of $T_f$ in Example \ref{example:biharmform-1}.}
        \end{figure}
    \end{center}
\end{example}

\begin{example}
    \label{example:biharmform-2}
    \par Consider the biharmonic form
    $$
        f(x, y) = - 448 x^5y^3 + 896 x^3 y^5 - 192 x y^7.
    $$
    Expressing $f$ in the form given in Observation \ref{observation:biharmform}, we obtain
    $$
        f = - 6i z^8 + 6i \bar{z}^8 + 8i z^7 \bar{z} - 8i z \bar{z}^7.
    $$
    Taking $\eta = 1$, we get $T_f = 6$. In this case $T_f$ is precisely a common cusp of $\Gamma_6$ and $\Pi_6$, and
    $$
        \res(P_{6, 6}, P^\#_{6, 6}) = \res(Q_{6, 6}, Q^\#_{6, 6}) = 0.
    $$
    Theorem \ref{theorem:maintheorem-dge4} therefore gives $\rank_\mathbb{R} f = 8$ and $\brank_\mathbb{R} f = 6$.

    \begin{center}
        \begin{figure}[htbp]
            \label{figure:example-tf-2}
            \begin{tikzpicture}[scale = 0.5]
    % Epicycloid E_6
    \path[fill = gray!25, draw = gray, line width = 1.0pt]
    plot[domain = 0:360, samples = 361, variable = \t]
    ({7*cos(\t) - cos(7*\t)}, {7*sin(\t) - sin(7*\t)}) -- cycle;

    % Hypocycloid H_6
    \path[fill = gray!50, draw = darkgray, line width = 1.0pt]
    plot[domain = 0:360, samples=361, variable=\t]
    ({5*cos(\t) + cos(5*\t)}, {5*sin(\t) - sin(5*\t)}) -- cycle;

    % Coordinate axes
    \draw[->, line width = 1.0pt] (-9,0) -- (9,0) node[right] {$\Re T$};
    \draw[->, line width = 1.0pt] (0,-9) -- (0,9) node[above] {$\Im T$};

    \foreach \x in {-8, -6, -4, -2, 2, 4, 6, 8}
    \draw (\x, 0.12) -- (\x, -0.12) node[below, font = \scriptsize] {$\x$};
    \foreach \y in {-8, -6, -4, -2, 2, 4, 6, 8}
    \draw (0.12, \y) -- (-0.12, \y) node[left, font = \scriptsize] {$\y$};

    % Tf
    \fill (6, 0) circle (4.5pt);
    \node[inner sep = 1.5pt] at (6, 0.8) {$T_f$};
\end{tikzpicture}
            \caption{The location of $T_f$ in Example \ref{example:biharmform-2}.}
        \end{figure}
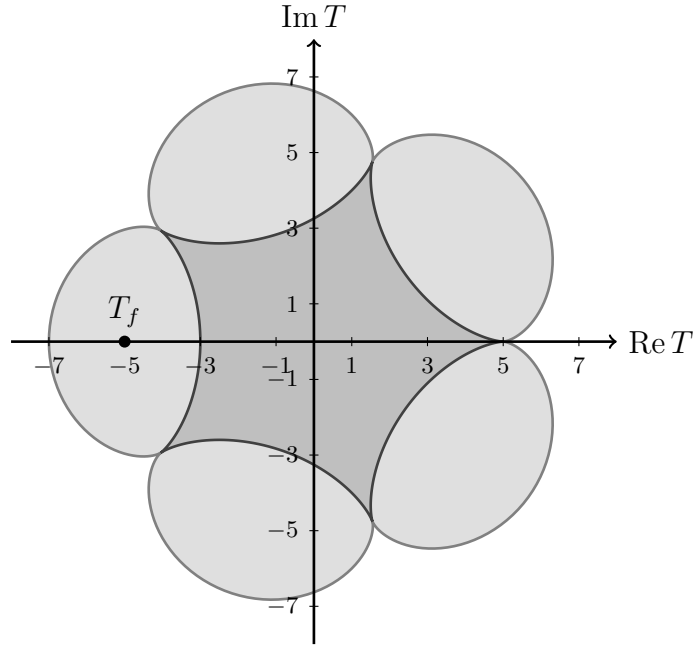
    \end{center}
\end{example}

    \bibliographystyle{abbrv}
    \bibliography{references}
\end{document}